%



\magnification=\magstep1
\input amstex
\documentstyle{amsppt}
\pagewidth{6.5truein}
\pageheight{8.9truein}
\NoRunningHeads

\def\C{\bold C}
\def\Ck{\bold C_{\kappa}}

\topmatter

\title
Semi-Cohen Boolean algebras 
\endtitle

\author
Bohuslav Balcar, Thomas Jech and Jind\v rich Zapletal
\endauthor

\affil
The Academy of Sciences of Czech Republic \\
The Pennsylvania State University
\endaffil

\thanks
Supported in part by a grant no. 11904 from AV\v CR (Balcar), by the National
Science Foundation grants INT-9016754 (Balcar and Jech,
U.S.--Czechoslovakia Cooperative Grant)
and DMS-9401275 (Jech and Zapletal) and by a National Research Council
COBASE grant (Jech). Balcar acknowledges the hospitality of the
Pennsylvania State University during his visit; Jech is grateful
for the hospitality of the Center for Theoretical Study in Prague\endthanks

\address
Mathematical Institute of the Academy of Sciences of Czech Republic,
\v Zitn\'a 25, Praha~1, Czech Republic (Balcar)
\endaddress
\email balcar\@earn.cvut.cz \endemail

\address
Department of Mathematics, The Pennsylvania State University,
University Park, PA 16802 (Jech, Zapletal)
\endaddress
\email jech\@math.psu.edu, zapletal\@math.psu.edu \endemail

\address
Center for Theoretical Study, Jilsk\'a 1, 110 00 Praha 1, Czech Republic
(current address of B. B. and T. J.)
\endaddress
\email jech\@ruk.cuni.cz \endemail

\address
Math. Sciences Research Inst., 1000 Centennial Drive, 
Berkeley, CA 94720 (current address of J.Z.)
\endaddress 

\abstract
We investigate classes of Boolean algebras related to the notion of forcing
that adds Cohen reals. A {\it Cohen algebra} is a Boolean algebra that is
dense in the completion of a free Boolean algebra. We introduce and study
generalizations of Cohen algebras: semi-Cohen algebras, pseudo-Cohen algebras
and potentially Cohen algebras. These classes of Boolean algebras are closed
under completion.
\endabstract

\endtopmatter

\document

\head{1.  Introduction} \endhead

For an infinite cardinal $\kappa$ let $\Ck$ denote the complete Boolean
algebra that adjoins $\kappa$ Cohen reals.  $\Ck$ is the completion of the
free Boolean algebra on $\kappa$ generators; equivalently, $\Ck$ is the
algebra of all regular open subsets of the topological product space
$\{0,1\}^{\kappa}$.  We call a Boolean algebra $B$ a {\it Cohen algebra}
if the completion of $B$ is $\Ck$.

We investigate Boolean algebras that closely resemble Cohen algebras,
particularly the class of algebras called {\it semi-Cohen}.
These algebras were introduced in 1992 by Fuchino and Jech, motivated
by Koppelberg's work on Cohen algebras \cite{Ko2}. The work
on this project was done between 1992 and 1995 during Balcar's visits
at Penn State and Jech's visits in Prague.

Semi-Cohen (called regularly filtered) and related Boolean algebras
are also the subject of a recent monograph by Heindorf and Shapiro
\cite{H-S}. Their work deals with generalizations of projectivity and
uses algebraic, rather than set theoretic, methods and point of view.

A game property equivalent to semi-Cohen
also appears in \cite{D-K-Z}.

An application of semi-Cohen algebras in topological dynamics appears in
\cite{B-F}.

\definition
{Definition 1.1}  A Boolean algebra $B$ of uniform density $\kappa$ is
{\it semi-Cohen} if $[B]^\omega$ has a closed unbounded set of countable
regular subalgebras of $B.$
\enddefinition

Every Cohen algebra is semi-Cohen.  In Section 3 we prove the following
characterization of Cohen algebras, a slight improvement of results
due to Koppelberg \cite{Ko2} and Bandlow \cite{Ba}:

\proclaim
{Theorem 3.2} \ Let $\kappa$ be an uncountable cardinal and let $B$ be a
Boolean algebra of uniform density $\kappa$.  Then $B$ is a Cohen algebra if
and only if the set
$$
	\{ A \in [B]^{\omega}:  A \le_{\text{reg}} B\}
$$
contains a closed unbounded set $C$ with the property
$$
	\text{if } A_1, A_2 \in C \text{ then }
	\langle A_1 \cup A_2 \rangle \in C
$$
where $\langle A_1 \cup A_2 \rangle$ is the subalgebra of $B$ generated by
$A_1 \cup A_2$.
\endproclaim

It is not difficult to see that every regular subalgebra, of uniform
density, of a semi-Cohen algebra is itself semi-Cohen (cf. Theorem
4.1).  In particular, every complete subalgebra of uniform density of
$\Ck$ is semi-Cohen.  Thus the concept of semi-Cohen algebras is
relevant to the following problems (cf. \cite{Ko1}, \cite{Ka}):

{\it Is every complete subalgebra of $\Ck$ of uniform density isomorphic
to some $\C_{\lambda}$?}

This problem was recently solved by Koppelberg and Shelah in
\cite{Ko-Sh} and we return to it below.

In Section 4 we investigate semi-Cohen algebras.  Among the results
of Section~4 are these (some have also been proved by Fuchino):

\proclaim
{Theorem 4.3} \ Let $B$ be a Boolean algebra of uniform density.
The following are equivalent:
\roster
\item"(a)"  $B$ is a semi-Cohen

\item"(b)"  $V^P \vDash$ $B$ is a Cohen, where $P$ is the $\sigma$-closed
collapse of $|B|$ onto $\aleph_1$.

\item"(c)"  There exists a proper forcing $P$ such that $V^P \vDash B$
is Cohen.

\item"(d)"  The second player has a winning strategy in the infinite game
$\Cal G$ in which two players select in turn elements
$a_0, b_0, a_1, b_1,\dots$ of $B$ and the second player wins if
$\langle \{a_0, b_0, a_1, b_1,\dots\} \rangle \;\le_{\text{reg}}$ $B$.

\item"(e)" \cite{D-K-Z} I has a winning strategy in the following game $\Cal H$:
I plays elements $a_i$ of $B$ and II plays elements $b_i\le a_i$. I wins
iff $\sum b_i = 1.$ 
\endroster
\endproclaim

\proclaim
{Theorem 4.5} \ Let $\{B_{\alpha}: \alpha < \vartheta\}$ be a
continuous increasing chain of semi-Cohen Boolean algebras such that
$B_{\alpha} \le_{\text{reg}} B_{\alpha+1}$ for every $\alpha$.  Then
$\underset{\alpha < \vartheta}\to \bigcup B_{\alpha}$ is semi-Cohen.
\endproclaim

\proclaim
{Theorem 4.6} \ (a) \ If $A$ and $B$ are semi-Cohen algebras then
$A \times B$ is  semi-Cohen.

(b) \ If $A$ is semi-Cohen and $V^A \vDash $ $\dot B$ is semi-Cohen then
$A * \dot B$ is semi-Cohen.

(c) \ If $B$ is semi-Cohen and if $W$ is an extension of the universe
then $W \vDash B$ is semi-Cohen.
\endproclaim

In Section 5 we present some examples:

\proclaim
{Theorem 5.2} \ There exists a semi-Cohen algebra of density
$\aleph_2$ that cannot be embedded as a regular subalgebra into a
Cohen algebra.
\endproclaim

\proclaim
{Theorem 5.11} \  There exists an increasing $\omega$-chain of
Cohen algebras $B_n$, with $B_n \le_{\text{reg}} B_{n+1}$ whose
union is not a Cohen algebra
\endproclaim

We also give a simplified proof of the result of Koppelberg and
Shelah mentioned above:

\proclaim
{Theorem 5.1} (Koppelberg-Shelah) \  For every $\kappa \ge \aleph_2$, $\Ck$
has a complete subalgebra of uniform density $\kappa$ that is not Cohen.
\endproclaim

We mention a related result from \cite{H-S}, 6.3.2: a construction of an
rc-filtered Boolean algebra (hence semi-Cohen) that is not a Cohen algebra.

\medpagebreak

In Section 6 we consider pseudo-Cohen algebras:

\definition
{Definition 6.1}  A Boolean algebra $B$ of uniform density $\kappa$ is
{\it pseudo-Cohen} if it has a stationary set of countable
regular subalgebras.
\enddefinition

Clearly, every semi-Cohen algebra is pseudo-Cohen, and every
regular subalgebra of a pseudo-Cohen algebra is pseudo-Cohen
(cf. Theorem 6.2).  We characterize pseudo-Cohen algebras, prove
preservation properties, and give some examples:

\proclaim
{Theorem 6.3} \ Let $B$ be an algebra of uniform density $\kappa$.  The
following are equivalent:

\roster
\item"(a)"  $B$ is a pseudo-Cohen.

\item"(b)"  There exists an $\aleph_0$-distributive forcing $P$ such that
$V^P \vDash B$ is Cohen.

\item"(c)"  The first player does not have a winning strategy in the
game $\Cal G$.
\endroster
\endproclaim

\proclaim
{Proposition 6.4} \ If $B$ is pseudo-Cohen and $W$ is a proper-forcing
extension of the universe then $W \vDash B$ is pseudo-Cohen.
\endproclaim

\proclaim
{Theorem 6.5} \ There exists a pseudo-Cohen algebra that is not
semi-Cohen.
\endproclaim

In Section 7 we consider a further generalization:

\definition
{Definition 7.1} A Boolean algebra $B$ of uniform density is
{\it potentially Cohen} if there exists a forcing $P$ preserving
$\aleph_1$ such that $V^P \vDash B$ is Cohen.
\enddefinition

By Theorem 6.2 (b), every pseudo-Cohen algebra is potentially Cohen.
As for the converse, we present two results:

\proclaim
{Theorem 7.2} \ The continuum hypothesis implies that every potentially
Cohen algebra is pseudo-Cohen.
\endproclaim

\proclaim
{Theorem  7.3} \ It is consistent (and follows from $MA + \neg CH$)
that the measure algebra (which is not pseudo-Cohen) is potentially
Cohen.
\endproclaim

\head{2. Preliminaries} \endhead

For a Boolean algebra $B$, we denote $B^+$ the set of all nonzero
elements.  We use +, $\cdot$ and - to denote Boolean-algebraic
operations and $\le$ for the Boolean algebraic ordering (inclusion).  Infinite
sums and products, when they exist, are denoted $\sum$ and $\prod$.  If
$B$ is a Boolean algebra and $X \subseteq B$, we denote
$$
	\langle X \rangle \;=
	\bigcap \{A:A \text{ is a subalgebra of $B$ and } X \subseteq A\}
$$
the subalgebra generated by $X$.  Every element of $\langle X \rangle$ can be
written as the sum $p_1+\cdots+p_n$ where each $p_i$ is
$\pm \,x_i \cdot \pm \,x_2 \cdot \dots \cdot \pm \,x_k$
with $x_1,\dots,x_k \in X$.  If $A$ is a subalgebra of $B$ and
$b_1,\dots, b_n \in B$,  we let
$$
	A (b_1,\dots,b_n) = \langle A \cup \{ b_1,\dots,b_n\}\rangle .
$$
Note that
$$
	A (b) = \{ a_1\cdot b + (a_2 - b): a_1, a_2 \in A\}.
$$

\definition
{Definition 2.1}  A subalgebra $A$ of $B$ is a
{\it regular subalgebra}
$$
	A \le_{\text{reg}} B,
$$
if for any $X \subseteq A$, if $\sum^A X$ exists then
$\sum^A X = \sum^B X$.
\enddefinition

The following equivalences are well known; cf. \cite{Ko1}:

\proclaim
{Lemma 2.2} \ The following are equivalent:

\roster
\item"(a)"  $A \le_{\text{reg}} B$,

\item"(b)"  every maximal antichain in $A$ is maximal in $B$,

\item"(c)"  for every $b \in B^+$ there exists some $a \in A^+$ such
that no $x \in A^+$ exists with the property that $x \le a$ and
$x \cdot b = 0$.
\endroster
\endproclaim

A set $D \subseteq B$ is {\it dense} in $B$ if for every
$b \in B^+$ there exists a $d \in D$ with $0 < d \le b$.  The
{\it density} of $B$ is the least size of a dense subset of $B$.
$B$ has {\it uniform density} if for every $a \in B^+$,
$B \upharpoonright a$ has
the same density.  (Every Boolean algebra can be decomposed into
algebras of uniform density.)  For every Boolean algebra $B$ there
exists a unique complete Boolean algebra $\overline B$, the
{\it completion} of $B$, such that $B$ is dense in $\overline B$.
The next lemma summarizes some known facts about regular subalgebras:

\proclaim
{Lemma 2.3}
\roster
\item"(a)"  If $A$ is a finite subalgebra of $B$ then
$A \le_{\text{reg}} B$.

\item"(b)"  If $A \le_{\text{reg}} B$ and $B \le_{\text{reg}} C$ then
$A \le_{\text{reg}} C$.

\item"(c)"  If $A$ is a subalgebra of $B$, $B$ is a subalgebra of $C$
and if $A \le_{\text{reg}} C$ then $A \le_{\text{reg}} B$.

\item"(d)"  If $A$ is a dense subalgebra of $B$ then
$A \le_{\text{reg}} B$.

\item"(e)"   $A \le_{\text{reg}} B$ if and only if
$\,\overline A \le_{\text{reg}} \overline B$.

\item"(f)"   If $A$ and $B$ are complete then
$A \le_{\text{reg}} B$ if and only if $A$ is a complete subalgebra
of $B$.

\item"(g)"   If $\{A_i\}_{i \in I}$ is a directed system of
algebras such that $A_i \le_{\text{reg}} A_j$ whenever $i \le j$,
and if $B = \underset{i \in I}\to\bigcup A_i$, then
$A_i \le_{\text{reg}} B$ for all  $i \in I$.
\endroster
\endproclaim

If $A$ is a subalgebra of $B$ and  $b \in B$, then the 
{\it (upper) projection} of $b$ to $A$ is the smallest element $a \in A$,
if it  exists, such that $b \le a$.  The projection of $b$ is
denoted $pr^A(b)$.

If $pr^A(b)$ exists for all $b \in A$, then
$A \le_{\text{reg}} B$, and the lower projection $pr_A(b)$ exists for
all $b \in B$, where  $pr_A(b) = $ the largest $a \le b$ in $A$.

A set $X \subseteq B$ is {\it independent} if
$$
	\pm \,x_1 \cdot \pm \,x_2 \cdot \dots \cdot \pm \,x_n \ne 0
$$
for all (distinct) $ x_1, \cdots, x_n \in X$.  
A Boolean algebra $B$ is {\it free}
over $X$ if $X$ is independent and $B = \langle X\rangle$.  The {\it free
algebra}
over $X$ is unique up to isomorphism and will be denoted $Fr_X$. We note
that $Fr_X$ is isomorphic to the set
algebra of all clopen sets of the Cantor space $\{0,1\}^X$.  We also
note that if $X \subseteq Y$ then $Fr_X \le_{\text{reg}} Fr_Y$.

If $A$ is an (uncountable) set, we denote $[A]^{\omega}$ the set of all
countable subsets of $A$, and $[A]^{ < \omega}$ the set of all finite
subsets of $A$.  $A$ set $C \subseteq [A]^{\omega}$ is {\it closed unbounded}
if $C$ is closed under unions of countable chains and for
every $x \in [A]^{\omega}$ there exists a $y \in C$ with $x \subseteq y$.
Some facts about closed unbounded sets:

\proclaim
{Lemma 2.4}

\roster
\item"(a)"  If $C \subseteq [A]^{\omega}$ is closed unbounded
then there exists a function $F: [A]^{<\omega}\to A$ such that
$$
	C \supseteq \{ x \in [A]^{\omega}: F (e) \in x \text{ whenever }
	e \in [x]^{<\omega}\}.
$$

\item"(b)"  If $C$ is closed unbounded and $D \subseteq C$  is
countable and directed under inclusion then $\bigcup D \in C$.

\item"(c)"  If $C \subseteq [A]^{\omega}$ is closed unbounded and
$A \subseteq B$, then $\{ x \in [B]^{\omega}: x \cap A \in C\}$ is closed
unbounded in $[B]^{\omega}$.

\item"(d)"  If $C$ is closed unbounded in $[B]^{\omega}$ and if
$A \subseteq B$ then $\{x \cap A : x \in C\}$ contains a closed
unbounded set in $[A]^{\omega}$.
\endroster

\endproclaim
If $(P, \le)$ is a notion of forcing then $B(P)$ will denote the
corresponding complete Boolean algebra, and $V^P = V^{B(P)}$ the
corresponding Boolean-valued model.

\head{3.  Cohen algebras} \endhead

For every infinite cardinal $\kappa$, let $\Ck$ be the complete Boolean
algebra that adjoins $\kappa$ Cohen reals.  In other words, $\Ck$ is the
completion of $Fr_{\kappa}$, the free algebra on $\kappa$
generators; more generally, let $\C_X = \overline{Fr_X}$ for any set $X.$

\definition
{Definition 3.1} \  A Boolean algebra $B$ is a {\it Cohen algebra}
if $\overline B = \Ck$ for some infinite cardinal $\kappa$.
\enddefinition

In this
section we prove the following characterization of Cohen algebras,
that is a slight improvement of results of Koppelberg \cite{Ko2} and
Bandlow \cite{Ba}.

\proclaim
{Theorem 3.2} \ Let  $B$ be an infinite Boolean algebra of uniform
density.  $B$ is a Cohen algebra if and only if the set
$\{A \in [B]^{\omega} : A\le_{\text{reg}} B \}$ contains a closed
unbounded set $C$ with the property
$$
	\text{if }  A_1, A_2 \in C, \text{ then }
	\langle A_1 \cup A_2 \rangle \;\in C.
\tag{$*$}
$$
\endproclaim

We remark that if $B$ is a countable atomless algebra 
then $B$ is a Cohen algebra, and
the condition is satisfied trivially, since $C= \{B\}$ is (trivially) a
closed unbounded subset of $[B]^{\omega}$.

We shall prove Theorem 3.2 in a sequence of lemmas. Throughout, we assume
that $B$ has uniform density.

\proclaim
{Lemma 3.3} \ If  $B$ is a dense subalgebra of $\Ck$ then
$B$ has  the property in Theorem 3.2.
\endproclaim

\demo
{Proof} \ Let $C$ be the set of all countable subalgebras $A$ of $B$
with the property that there exists a countable set $S \subseteq \kappa$
such that $A$ is dense in $B \cap \C_S$ and $B \cap \C_S$ is dense in
$\C_S$.  We shall prove that every $A \in C$ is a regular subalgebra of
$B$, that $C$ is closed unbounded and that ($*$) is satisfied.

Let $A \in C$, with $S$ being a witness.  Since $B \cap \C_S$ is
dense in $\C_S$ and $\C_S \le_{\text{reg}} \Ck$, we have
$B \cap \C_S \le_{\text{reg}} \Ck$, and because $B$ is dense in
$\Ck$ it follows that $B \cap \C_S \le_{\text{reg}} B$.  Since $A$
is dense in $B \cap \C_S$, we have $A \le_{\text{reg}} B$.

To show that $C$ is unbounded, let $a \in B$ be arbitrary; we shall
find an $A \in C$ such that $a \in A$.  First, because $\Ck$ has
the countable chain condition, there exists a countable $S_0 \subseteq
\kappa$ such that $a \in \C_{S_0}$.  Second, again using the countable
chain  condition, and because $B$ is dense in $\Ck$, there exists a
countable $S \subseteq \kappa$ such that $S_0 \subseteq S$ and that
$B \cap \C_S$ is dense in $\C_S$.  Finally, there is a countable subalgebra
$A$ of $B$ such that $a \in A$ and that $A$ is dense in $B \cap \C_S$.

To show that $C$ is closed, let $\{A_n\}^{^\infty}_{_{n=0}}$
be an increasing chain
in $C$ and let $A = \overset\infty\to{\underset{n=0}\to \bigcup} A_n$.  Let
$\{S_n\}^{^\infty}_{_{n=0}}$ be witnesses for the $A_n$.
For each $n$, $A_n$ is
dense in $\C_{S_n}$ and a subalgebra of $\C_{S_{n+1}}$; hence
$S_n \subseteq S_{n+1}$.  Since for each $n$, $A_n$ is dense in $\C_{S_n}$,
it follows that $A$ is dense in
$\overset\infty\to{\underset{n=0}\to \bigcup} \C_{S_n}$;
but the latter is
dense in $\C_S$ where $S = \overset\infty\to{\underset{n=0}\to \bigcup} S_n$.
Hence $A$ is dense in $B \cap \C_S$ which is dense in $\C_S$.

Finally, we shall verify ($*$).  Let $A_1, A_2 \in C$ and let $S_1, S_2$
be such that $A_i$ is dense in $\C_{S_i}$, for $i = 1,2$.  We shall show
that $A = \;\langle A_1 \cup A_2 \rangle$ is
dense in $\C_S$ where $S = S_1 \cup S_2$.  Let
$b \in \C^+_S$ be arbitrary; we shall find $a_1 \in A_1$ and $a_2 \in A_2$
such that $0 \ne a_1 \cdot a_2 \le b$.  

The algebra $\C_S$ has as a dense set Cohen's forcing $P_S,$ the set of all
finite 0-1-functions on $S.$ Let $p\in P_S$ be such that $p\le b.$ Let
$p_1 = p \restriction S_1$ and $p_2= p\restriction S_2.$ First we find some
$a_1 \in A_1^+$ such that $a\le p_1$, and then some $q_1 \in P_{S_1}$ such
that $q_1\le a_1.$ Let $q_2 = p_2 \cup (q_1 \restriction S_2);$ we have
$q_2 \in P_{S_2}.$ Now we find some $a_2 \in A_2^+$ such that $a_2 \le q_2.$
We claim that $a_1\cdot a_2 \ne 0:$ there exists some $r_2\in P_{S_2}$ such
that $r_2\le a_2,$ and then $r_2 \cup q_1 \in P_S$ is below both $a_1$ and
$a_2.$   $\square$
\enddemo

This lemma proves one direction of Theorem 3.2.  For the other
direction, we first prove that the property in Theorem 3.2
implies that $B$ has the countable chain condition.  In fact,
we prove a stronger assertion and use it in Section 4.

\proclaim{Lemma 3.4} \  Let $B$ be a Boolean algebra such
$\{A \in [B]^{\omega} : A\le_{\text{reg}} B \}$ is stationary.
Then $B$ has the countable chain condition.
\endproclaim

\demo
{Proof} \ Let $W$ be a maximal antichain in $B$. Consider the model
$M=\langle  B, \le, W \rangle$.  There exists a stationary set of
countable submodels $A \prec M$ such that
$A\le_{\text{reg}} B$.  It follows that $W\cap A$ is a
maximal antichain in $A$ and therefore in $B$.  Hence
$W = W\cap A$ and is countable.   $\square$
\enddemo

The following lemma is due to Vladimirov \cite{V}, Lemma VII.3.

\proclaim
{Lemma 3.5} \  Let $B$ be a complete Boolean algebra
of uniform density and let $A$ be a complete subalgebra of $B$ of
density less than the density of $B$.  
For every
$v \in B$ there exists a $u \in B$ such that
$a \cdot u \ne 0 \ne a - u$ for all $a \in A^+$, and
$v \in A(u)$.
\endproclaim

We shall call $u$ {\it independent over} $A$.

\proclaim
{Lemma 3.6}  Let $B$ have a closed unbounded set $C$ of countable
regular subalgebras, closed under $\langle A_1 \cup A_2 \rangle$.  Let
$\Cal S$ be a collection of all subalgebras of $B$ of the form
$\langle \bigcup S \rangle$
where $S \subseteq C$.  Then every $A \in \Cal S$ is a
regular subalgebra of $B$.
\endproclaim

\demo
{Proof} \ This is true for every finite $S \subset C$, and because
$C$ is closed, it is also true for every countable $S \subset C$.
For an arbitrary $S$, let $W$ be a maximal antichain in
$\langle \bigcup S \rangle$.  By Lemma 3.4, $W$ is countable, and so
$W \subset \langle \bigcup S_0 \rangle$ for some countable $S_0 \subseteq S$.
Since $\langle \bigcup S_0 \rangle \;
\le_{\text{reg}} B$, $W$ is a maximal antichain
in $B$.  Hence $\langle \bigcup S \rangle \;\le_{\text{reg}} B$.  $\square$
\enddemo

We shall now complete the proof of Theorem 3.2.  Let $B$ be an
infinite Boolean algebra of uniform density $\kappa > \omega$ and let $C$
be a closed unbounded subset of $[B]^{\omega}$ consisting of regular
subalgebras of $B$ and satisfying ($*$).  Let
$$
	\Cal S = \{ \langle \bigcup S \rangle : S \subseteq C \},
$$
and let $D = \{ d_{\alpha} : \alpha < \kappa \} \subset B^+$ be a dense
subset of $B$.

If $A_1$ and $A_2$ are subalgebras of $\overline B$ we say that
$A_1$ and $A_2$ are {\it co-dense} if for every $a_1 \in A^+_1$
there exists some $a_2 \in A^+_2$ with $a_2 \le a_1$, and for
every $a_2 \in A^+_2$ there exists some $a_1 \in A^+_1$ with
$a_1 \le a_2$.

We shall construct, by induction on $\alpha < \kappa$, two continuous
chains $G_0 \subset G_1
\subset \dots \subset G_{\alpha}
\subset \dots$ and $B_0 \subset B_1 \subset \dots \subset
B_{\alpha} \subset \dots$ such that
\roster
\item  $B_{\alpha} \in \Cal S$,

\item  $G_{\alpha}$ is an independent subset of $\overline B$,

\item  $A_{\alpha} = \langle G_{\alpha} \rangle$ and $B_{\alpha}$ are co-dense,

\item  $d_{\alpha} \in B_{\alpha+1}$,

\item  $G_{\alpha+1} - G_{\alpha}$ is countable.
\endroster

This will prove that $B$ is a Cohen algebra, because by (4),
$\bigcup_{\alpha < \kappa} B_{\alpha}$ is dense in $B$,
hence $\bigcup_{\alpha < \kappa} A_{\alpha}$ is dense
in $\overline B$, and it follows that $\overline B =
\overline{Fr_G}$ where $G = \bigcup_{\alpha < \kappa}
G_{\alpha}$.

At limit stages of the construction, we let $B_{\alpha} =
\bigcup_{\beta < \alpha} B_{\beta}$ and
$G_{\alpha} = \bigcup_{\beta <\alpha} G_{\beta}$;
clearly, (1), (2) and (3) are satisfied.

Thus assume that we have constructed $B_{\alpha}$ and
$G_{\alpha}$, and find $G_{\alpha+1}$ and $B_{\alpha+1}$.  Since
$A_{\alpha}$ is dense in $\overline{B_{\alpha}}$  and
$\overline{B_{\alpha}}$ is a complete subalgebra of
$\overline B,$  $\overline A_{\alpha}$ is a complete subalgebra of
$\overline B$.  Also, if $u_1,\dots,u_n \in \overline B$ then
$\overline A_{\alpha} (u_1,\dots,u_n)$ is a complete algebra of
$\overline B$.

Since $|A_{\alpha}|< \kappa$ and $\kappa$
is the uniform density of $B$, we find, by Lemma 3.5,
for every $b \in \overline B$ some $u \in \overline B$ independent
over $\overline A_{\alpha}$ such that
$b \in \overline A_{\alpha}(u)$.  More generally, if $b,
u_1,\dots,u_n \in \overline B$ then there exists some $u$
independent over $\overline A_{\alpha}(u_1,\dots,u_n)$ such that
$b \in \overline A_{\alpha} (u_1,\dots,u_n, u)$.

Given $u \in \overline B$, there exist countably many
$\{b_n\}^{^\infty}_{_{n=0}} \subset B$ such that
$\overset\infty\to{\underset{n=0}\to\sum} b_n = u$.  Then there
exists some $X \in C$ such that $\{b_n\}_n \subseteq X$ and so
$\langle B_{\alpha} \cup X \rangle$ is dense in $\overline A_{\alpha}(u)$.
Therefore, there exist a countable set $\{u_n\}^{^\infty}_{_{n=0}}
\subset \overline B$, and some $B_{\alpha + 1} \in \Cal S$ such
that $d_{\alpha} \in B_{\alpha+1}$, that $G_{\alpha+1} =
G_{\alpha} \cup \{u_n\}^{^\infty}_{_{n=0}}$ is independent and that
$A_{\alpha+1} = \langle G_{\alpha+1} \rangle$
and $B_{\alpha + 1}$ are co-dense.

\head{4.  Semi-Cohen algebras} \endhead

Motivated by the characterization of Cohen algebras, Fuchino and
Jech introduced in 1992 the following property:  a Boolean
algebra $B$ of uniform density is called semi-Cohen if it has a closed
unbounded set of countable regular subalgebras.

We start with the following observation:

\proclaim
{Theorem 4.1} \ If $B$ is semi-Cohen and if $A$ is a regular
subalgebra of $B$ of uniform density then $A$ is semi-Cohen.
\endproclaim

\demo
{Proof} The family $[B]^\omega$ contains a closed unbounded subset of 
regular subalgebras
of $B$.  Since $A\le_{\text{reg}} B$, there exists for every
$b \in B^+$ some $a \in A^+$ such that $a-b$ does not have any
$x \le a-b$ in $A^+$.  Let $F:B^+ \to A^+$ be a function that to
each $b \in B^+$ assigns such an $a \in A^+$.  Let $C \subseteq
[B]^{\omega}$ be a closed unbounded set of regular subalgebras, closed under
$F$.

If $X \in C$ then $A \cap X \le_{\text{reg}} X$ because $X$ is
closed under $F$.  Every maximal antichain in $A\cap X$ is
maximal in $X$, hence in $B$ (because $X \le_{\text{reg}} B$), hence
in $A$.  So $A \cap X \le_{\text{reg}} A$.

There is a closed unbounded set $D \subseteq [A]^{\omega}$ such that
$D \subseteq \{ X \cap A : X \in C\}$.  $D$ witnesses that $A$ is
semi-Cohen.  $\square$
\enddemo

Another consequence of Definition 1.1 is

\proclaim
{Theorem 4.2} \ If $B$ is a semi-Cohen algebra and $B$ has density
$\aleph_1$ then $B$ is a Cohen algebra.
\endproclaim

\demo
{Proof} \ If $B$ is semi-Cohen, and if $|B| = \aleph_1$ then
$B$ is Cohen, because every closed unbounded subset of $[B]^{\omega}$
contains a closed unbounded subset that is a chain under inclusion,
hence closed under finite unions.  If $B$ has density $\aleph_1$,
let $A$ be a dense subalgebra of $B$ of size $\aleph_1$.  $A$ is
semi-Cohen by Theorem 4.1 and therefore Cohen, and since
$\overline A = \overline B$, $B$ is Cohen.  $\square$
\enddemo                        

As a consequence, we have the following corollary: If $B$ is a Cohen algebra 
and $A\le_{\text{reg}} B$
has uniform density $\aleph_1$, then $A$ is a Cohen algebra. This fact was
previously known to Koppelberg \cite{Ko}.

We shall now give equivalent characterizations of semi-Cohen
algebras.

\proclaim
{Theorem 4.3} \ Let $B$ be a Boolean algebra of uniform density.  The
following are equivalent:
\roster
\item"(a)"  $B$ is semi-Cohen.

\item"(b)"  $V^P \vDash B$ is Cohen, where $P$ is the $\sigma$-closed
collapse of $|B|$ onto $\aleph _1$.

\item"(c)"  There exists a proper forcing $P$ such that $V^P \vDash
B$ is Cohen.

\item"(d)"  The second player has a winning strategy in the infinite
game $\Cal G$ in which two players select in turn elements
$a_0, b_0, a_1, b_1,\dots$ of the algebra $B$, and the second player wins if
$\langle \{a_0, b_0, a_1, b_1,\dots\} \rangle \;\le_{\text{reg}} B$.

\item"(e)" I has a winning strategy in the following game $\Cal H:$ I plays
elements $a_i$ of $B$ and II plays elements $b_i\le a_i.$ I wins iff
$\sum b_i = 1.$

\endroster
\endproclaim

The game $\Cal H$ is introduced in \cite{D-K-Z} where the equivalence of
(e) with semi-Cohen is proved.

First we state a corollary of this theorem:

\proclaim
{Corollary 4.4} An algebra $B$ is semi-Cohen if and only if $\overline B$ is
semi-Cohen.
\endproclaim

\demo
{Proof} \ One direction follows from Theorem 4.1, since
$B\le_{\text{reg}} \overline B$.  Thus assume that $B$ is semi-Cohen,
and show that $\overline B$ is.

Let $A = \overline B$. Let $P$ be the collapse (with countable conditions)
of $|A|$ onto
$\aleph_1$.  In $V^P$, $B$ is dense in $A$ and is a Cohen algebra;
hence $A$ is Cohen in $V^P$.  Hence $A$ is semi-Cohen in $V$.  $\square$.

\medpagebreak
To prove Theorem 4.3 we first prove that (a) and (b) are equivalent.
Let $P$ be the collapse of $|B|$ onto $\aleph_1$ with countable
conditions.  By Theorem 4.2, if $B$ is semi-Cohen in $V^P$ then $B$ is
Cohen.  Thus it suffices to show that $V^P \vDash$``$ B$ is semi-Cohen'' if
and only if $B$ is semi-Cohen.

As $P$ does not add new countable sets of ordinals, $[B]^{\omega}$ is the
same in $V^P$ as in $V$.  Using property (c) of Lemma 2.2 we see that
for every $A \in [B]^{\omega}$, $A\le_{\text{reg}} B$ if and only if $V^P
\vDash A\le_{\text{reg}} B$.

Let $S$ be the set of all countable regular subalgebras of $B$.  If
$S$ contains a closed unbounded set $C$, then $V^P \vDash C$ is
closed unbounded.  Conversely, if $S$ does not contain a
closed unbounded set, then because $P$ is proper, $V^P \vDash S$
does not contain a closed unbounded set.  Hence $B$ is semi-Cohen
if and only if $V^P \vDash B$ is semi-Cohen.

A similar argument establishes the equivalence of (a) with (c):  As
(c) follows from (b), let us assume that $P$ is proper and
$V^P \vDash B$ is semi-Cohen.  Let $S$ be the set of all countable
regular subalgebras.  If $B$ were not semi-Cohen then $[B]^{\omega}-S$
would be stationary, therefore stationary in $V^P$, contrary to the
assumption that $B$ is semi-Cohen in $V^P$.

To see that (a) and (d) are equivalent, consider the
game $\Cal G$.  If $F:~B^{<\omega} \to B$ is a winning strategy for the
second player then the set $C$ of all countable subalgebras of $B$
closed under $F$ is closed unbounded and all of its elements are
regular subalgebras.  Conversely, if $[B]^{\omega}$ has a closed unbounded
set of regular subalgebras then it is easy to find a winning
strategy for the second player, using Lemma 2.4(a).

For the convenience of the reader we outline the proof of the equivalence
of (a) and (e), which is the content of Theorem 1.6 of \cite{D-K-Z}.
If I has a winning strategy $\sigma$ in $\Cal H,$ then the club $C$ of
all countable subalgebras of $B$ closed under $\sigma$ consists of regular 
subalgebras: if $A\in C$ were not regular then there would be some 
maximal antichain $W$ in $A$ and some $b\in B^+$ incompatible with $W.$
Then II can play moves $b_i\le a_i$ within $A$ incompatible with $b.$
In the end, $\sum b_i \perp b,$ contradicting the assumption on $\sigma.$

Conversely, if there is a club $C$ consisting of regular subalgebras of $B$,
I wins the game by catching her tail: make sure that $A=\{a_i\}_i =
\{b_i\}_i \in C.$ Then in $A$, $\sum b_i =1$, and by regularity,
$\sum b_i = 1$ in $B.$
$\square$
\enddemo

Next we shall prove that the class of semi-Cohen algebras is
closed under unions of regular chains.  In the next Section we
show that this is not necessarily true for Cohen algebras.

\proclaim
{Theorem 4.5} \ Let $\{B_{\alpha}: \alpha < \vartheta\}$ be a
continuous increasing chain of semi-Cohen Boolean algebras such
that $B_{\alpha} \le_{\text{reg}} B_{\alpha + 1}$ for every
$\alpha$.  Then $\underset{\alpha < \vartheta}\to\bigcup B_{\alpha}$
is semi-Cohen, provided it has uniform density.
\endproclaim

Let us remark that uniform density need not be preserved by limits of chains.
The theorem is analogous to
\v S\v cepin's Theorem on openly generated Boolean algebras \cite{Fu}.

\demo
{Proof} \ Let $\{B_{\alpha}: \alpha < \vartheta\}$ be a regular continuous
chain of semi-Cohen algebras with limit $B=B_\vartheta.$
By Lemma 2.2 we can choose
{\it pseudo-projections} $\pi_{\alpha \beta} : B_\beta \to B_\alpha$
for all $\alpha \le \beta \le \vartheta$, i.e. functions 
such that for every $b\in B_\beta^+$
and for every $x\in B_\alpha^+$, if $x\le \pi_{\alpha \beta}(b)$ then
$x\cdot b\ne 0$.

Let $\lambda$ be a sufficiently large regular cardinal.
We will show that for every countable elementary 
submodel $M$ of $\langle H_\lambda, \in,
B, \{B_\alpha\}_\alpha, \{\pi_{\alpha \beta}\}_{\alpha \beta}
\rangle$, $B\cap M$ 
is a regular subalgebra of $B.$ This will show that $B$ is semi-Cohen.

Fix the model $M$ and let $\gamma=\sup M \cap \vartheta.$ For every $b\in B^+$
we must find an $a\in B\cap M$ 
such that every nonzero $x\in B\cap M$ below $a$ is 
compatible with $b.$ Fix some 
$b\in B^+.$ Let $b_\gamma=\pi_{\gamma \vartheta}(b);$
there is some $\alpha\in\gamma\cap M$ such that $b_\gamma\in B_\alpha.$
Since $B_\alpha$ is semi-Cohen, $B_\alpha \cap M$ is a regular subalgebra
of $B_\alpha$ and there is an $a\in B_\alpha\cap M$ such that every nonzero
$z\in B_\alpha\cap M$ below $a$ is compatible with $b_\gamma$. We claim 
that this $a$ works.

To show that, let $x\in B\cap M$ 
be below $a$; we will prove that $x\cdot b\ne 0.$
Let $y=\pi_{\alpha \vartheta}(x);$ $a\cdot y$ is a nonzero element of
$B_\alpha\cap M$, $a\cdot y \le a$ and so $a\cdot y\cdot b_\gamma \ne 0$
by the choice of $a.$ Since $y=\pi_{\alpha \vartheta}(x)$ and $a\cdot y\cdot
b_\gamma \le y,$ we conclude that $a\cdot y\cdot b_\gamma \cdot x \ne 0.$
By the same reasoning $a\cdot y\cdot b_\gamma \cdot x\cdot b \ne 0$ and
therefore $x\cdot b \ne 0$, as needed.
$\square$

\enddemo

The next theorem lists other closure properties of the class of
semi-Cohen algebras.

\proclaim
{Theorem 4.6} \ (a) \ If $A$ and $B$ are semi-Cohen algebras then
$A \times B$ is semi-Cohen.

(b) \ If $A$ is semi-Cohen and $V^A \vDash \dot B$ is semi-Cohen
then $A * \dot B$ is semi-Cohen.

(c) \ If $B$ is semi-Cohen and $W$ is an extension of the universe
then $W \vDash B$ is semi-Cohen. \endproclaim

By (b) and by Theorem 4.5, any
finite support iteration of semi-Cohen algebras is semi-Cohen, provided
it has uniform density.
Note also that the class of all Cohen algebras also has properties
(a), (b) and (c).

\demo
{Proof} \ (a) is an easy consequence of property (b) of
Theorem 4.3.  A similar argument proves (b):  Let $P$ be the
$\sigma$-closed collapse of $|A * \dot B|$, and note that $\dot B$
is Cohen in $(V^P)^A$; thus $A * \dot B$ is Cohen in $V^P$.

To prove (c), we use property (e) of Theorem 4.3.  Let $\sigma$ be a
winning strategy for the first player in the game $\Cal H$.  The
following statement is a theorem of ZFC:

$$
\aligned
&\sigma \text{ is a winning strategy if and only if for every $b\in B$ 
the tree}\\
&T_b = \{g : g \text{ is a finite play of the game 
$\Cal H$ according to } \sigma\\
&\text{such that II's answers are incompatible with } b\} \\
&\text{is well-founded.}
\endaligned
\tag{*}
$$ 

For example, if $\sigma$ is not  a winning strategy then there is a play
$a_0, b_0,a_1, b_1,...$ of the game $\Cal H$ according to $\sigma$ such that
$\sum b_i \ne 1.$ Consequently, if $b$ is incompatible with
$\sum b_i$ then this play is an infinite path through the tree $T_b.$

Now the statement (*) is absolute between transitive models of set
theory. Therefore any winning strategy $\sigma$ for I remains a
winning strategy in any extension $W$ of the universe and so
$W \vDash B $ is semi-Cohen. $\square$
\enddemo

The forcing properties of semi-Cohen algebras are much like those of Cohen
algebras; e.g. only Cohen reals are added. We include the following lemma
that will be needed in Section 5:

\proclaim
{Lemma 4.7} Semi-Cohen forcings do not add new branches to trees of
height $\omega_1.$
\endproclaim

\demo
{Proof} \ For contradiction, let $P$ be semi-Cohen, let $T$ be a tree
of height $\omega_1$ and let $p\in P$ force that $\dot b$ is a new branch
through $T.$ Choose a sufficiently large regular cardinal $\lambda$
and a countable elementary submodel $M$ of $H_\lambda$ that contains
$P, p, \dot b, T.$ Let $\gamma=M\cap\omega_1$ and choose a condition
$q\le p$ and some $t$ at level $\gamma$ of $T$ such that $q \Vdash
t\in \dot b.$ Since $P$ is semi-Cohen, the poset $P\cap M$ is regular
in $P$ and there is an $r\in P\cap M$ such that any extension of $r$
in $M$ is compatible with $q.$ But then, $c=\{s\in T : \text{ there is
some extension of $r$  forcing $s$ into } \dot b\}$ is a branch through $T$
and so $r\Vdash \dot b = c$, contradicting our assumption. $\square$
\enddemo

It should be remarked that the lemma can be easily adapted for the
two generalizations of semi-Cohenness in Sections 6 and 7.

\medpagebreak

The last result of this Section is the following theorem:

\proclaim
{Theorem 4.8} \ Every semi-Cohen complete Boolean algebra
$B$ of density $\kappa$ contains $\Ck$ as a complete subalgebra.
\endproclaim

\demo
{Proof} \ Let $B$ be a semi-Cohen complete Boolean algebra of uncountable
density $\kappa$. (The case when $\kappa= \aleph_0$ is trivial.) 
There exists a function $F: B^{<{\omega}}\to B$
such that every countable subalgebra of $B$ closed under $F$ is a
regular subalgebra.  First we claim that {\it every} subalgebra
of $B$ closed under $F$ is a regular subalgebra.  Let $A$ be a
subalgebra of $B$ closed under $F$ and let $W$ be a maximal
antichain in $A$.  Let $A_1$ be the smallest subalgebra of $B$
closed under $F$ such that $W \subset A_1$. Then  $A_1$ is a
subalgebra of $A$, and $W$ is a maximal antichain in $A_1$.  As
$A_1$ is countable, it is a regular subalgebra of $B$ and so $W$
is a maximal antichain in $B$.

Using $F$, we can find a continuous chain of regular subalgebras of
$B$ of size less than $\kappa$, whose union is dense in $B$.  Thus
there exist complete Boolean subalgebras  $B_{\alpha}$ of $B$,
$\alpha < \kappa$, each of density $< \kappa$, such that $B_{\alpha}\subset
B_{\beta}$ whenever $\alpha < \beta$, and for
every limit $\lambda$,
$\underset{\alpha < \lambda}\to\bigcup B_{\alpha}$ is dense in
$B_{\lambda}$.  As $B$ has uniform density $\kappa$, we may assume that
for every $\alpha$, $V^{B_{\alpha}}\vDash B_{\alpha+1}:B_{\alpha}$
is nontrivial.  Thus $B = B(P)$ where $P$ is the finite support
iteration of $\langle \dot Q_{\alpha} : \alpha < \kappa \rangle$, with
$\dot Q_{\alpha} = B_{\alpha+1}:B_{\alpha}$.

We now claim that $P$ (in fact every finite support iteration
of nontrivial forcings) embeds $\Ck$.  For each $\alpha$, let
$\dot a_{\alpha} \in V^{B_{\alpha}}$ be a name for an element of
$B_{\alpha+1}:B_{\alpha}$ such that
$V^{B_{\alpha}}\vDash 0 \ne \dot a_{\alpha} \ne 1$.

If $p \in P$ then we can find some stronger $q$ with the property
that for every $\alpha \in \text{support}(q)$, $q$ decides
$\dot a_{\alpha} \in \dot G$, where $\dot G$ is the name for the
generic ultrafilter.  Let $Q$ be the dense subset
of $P$ consisting of
all $q$ with such a property.  For every $q \in Q$, let $h (q)$ be
the 0-1 function on support$(q)$ such that $h(q) (\alpha)=1$
just in case $q \Vdash \dot a_{\alpha} \in \dot G$.

Let $P_{\kappa}$ be the (version of) Cohen forcing consisting of all
finite 0-1 functions $f$ with dom$(f) \subset \kappa$; we have
$B(P_{\kappa})=\Ck$.  The function $h$ maps $Q$ onto $P_{\kappa}$ and has the
property that for every $q \in Q$ and every $f \supset h(q)$ there
is some $q' \le q$ such that $h(q') \supseteq f$.  Hence $h$ witnesses
that $P_{\kappa}$ embeds regularly into $Q$, and so $\Ck$ embeds into
$B$ as a complete subalgebra.  $\square$
\enddemo

\head{5. Examples} \endhead

By Theorem 4.2, there is only one complete semi-Cohen algebra
of uniform density $\aleph_1,$ namely the Cohen algebra $\C_{\omega_1}.$
In this section we show that for higher densities
there are other, significantly different, semi-Cohen algebras.

Since semi-Cohenness is inherited by complete subalgebras,
we can find semi-Cohen algebras
looking at complete subalgebras of $\C_{\omega_2}$.  But are not all
of these again isomorphic to $\C_{\omega_2}$? (\cite{Ka}, \cite{Ko1}.)  The
following result from \cite{Ko-Sh} shows that this is not the case
and thereby provides a new kind of semi-Cohen algebras.

\proclaim
{Theorem 5.1} (Koppelberg-Shelah)
There is a complete subalgebra of $\C_{\omega_2}$ of
uniform density $\aleph_2$ which is not isomorphic to $\C_{\omega_2}$.
\endproclaim

\demo
{Proof} \ Let us define a partially ordered set $P$ as follows:
$$
	P = \{ z: z \text{ is a function, dom}(z)
	\in [\omega_2]^{<\omega}, \text{ ran}
	(z)\subseteq \omega^{< \omega}\}.
$$
We order $P$ by $z \le w$ if $z$ is a coordinatewise extension of
$w$ and for $\alpha \ne \beta$ both in dom$(w)$ if $n \in
\text{dom}(z(\alpha) - w(\alpha))$ and $n \in
\text{dom}(z(\beta))$ then $z(\alpha)(n) \ne  z(\beta)(n)$.  Thus
$P$ is the forcing for adding a sequence of $\omega_2$
eventually different reals.  If $G \subseteq P$ is a generic filter
then $G$ can be decoded from a $P$-generic sequence of functions
$\langle f_{\alpha}: \alpha < \omega_2 \rangle
\subseteq \omega^\omega$,
where $f_{\alpha} = \bigcup\{ z(\alpha): z\in G \}$.

Let $B =B(P)$.  We shall prove that the algebra $B$
witnesses  the statement of the theorem.  First, a helpful
observation.
\enddemo

\proclaim
{Lemma 5.2} \ The algebra $B$ is isomorphic to the completion of
the poset $R$, a finite support iteration
$$
	R = \;\langle R_{\alpha}: \alpha \le \omega_2, \dot Q_{\alpha}:
	\alpha < \omega_2 \rangle
$$
such that for every  $\alpha < \omega_2 $ the poset $R_{\alpha}$
forces:

\roster
\item  The $Q_{\alpha}$ generic is given by a real $\dot f_{\alpha}$,

\item  $Q_{\alpha} = \{ \langle s, a \rangle \; : s \in \omega^{<\omega}$,
$a \in [\alpha]^{< \omega}\}$ ordered by
$\langle t, b \rangle \, \le \,\langle s, a \rangle$ if
$s \subseteq t$, $a \subseteq b$ and $\forall  n\in \text{dom }(t - s)$
$\forall \; \beta \in a \quad t(n) \ne f_{\beta} (n)$,

\item  the $Q_{\alpha}$-term  $\dot f_{\alpha}$ is defined by
$Q_{\alpha} \Vdash$ ``$\dot f_{\alpha} = \bigcup \{ s \in \omega^{<\omega}:
\; \langle s, \emptyset \rangle \, \in \dot G$\}'', where $\dot G$ is the name for the
$Q_{\alpha}$-generic filter,

\item  $Q_{\alpha}$ is a separative partial order of uniform density
$|\alpha| + \aleph_0$.
\endroster
\endproclaim

\demo
{Proof} Left to the reader.  $\square$
\enddemo

Following the Lemma, we shall represent the
algebra $B$ as $B = \bigcup_{\alpha < \omega_2} B_{\alpha}$, where
$B_{\alpha} = B(R_{\alpha})$.  The following two lemmas
finish the proof of Theorem 5.1.

\proclaim
{Lemma 5.3} \  $B$ is not isomorphic to $\C_{\omega_2}$.
\endproclaim

\proclaim
{Lemma 5.4} \  $B$ can be completely embedded into $\C_{\omega_2}$.
\endproclaim

\demo
{Proof of Lemma 5.3} \ For contradiction, assume that
$h: B\to \C_{\omega_2}$ is an isomorphism.  Then by a simple
closure argument, there is an $\alpha,$ $\omega_1 < \alpha < \omega_2$,
such that $h'' B_{\alpha} = \C_{\alpha}$.
We reach a contradiction working in $V^{\C_{\alpha}}$. Let
$Q_{\alpha}, \dot f_{\alpha}$ be as defined in Lemma 5.2; we have
$Q_{\alpha} \in V^{\C_{\alpha}}$.  Since the residue forcing
$\C_{\omega_2} : \C_{\alpha}$ is Cohen, every real added by it comes
from a $\C_\omega$-extension of $V^{\C_{\alpha}}$.  In particular, the
real $\dot f_{\alpha}$ comes from such an extension.  Since
$\dot f_{\alpha}$ determines a $Q_{\alpha}$-generic filter
(Lemma 5.2 (1)) we have $B(Q_{\alpha}) \le_{\text{reg}} \C_\omega$.
This is a contradiction, since  $B(Q_{\alpha})$ has uniform
density $\aleph_1$ (Lemma 5.2 (4)) and  $\C_\omega$ has uniform density
$\aleph_0$.   $\square$
\enddemo

\demo
{Proof of Lemma 5.4} \ We define another forcing $S$ by

$S = \{ w: w \text{ is a finite function with dom}(w) \subseteq
(\omega_2 \times 2) \cup \{\omega_2\}$ such that
\roster
\item {}{}\quad $\langle \alpha, 0 \rangle \in \text{dom}(w)
\Rightarrow w(\alpha, 0) \in \omega^{<\omega},$

\item {}{}\quad $\langle \alpha, 1 \rangle \in \text{dom}(w)
\Rightarrow  w(\alpha, 1) \in \omega$,

\item {}{} \quad $ \omega_2 \in \text{dom}(w) \Rightarrow
w(\omega_2)$ is a finite function from  $\omega^{<\omega}$  which is
 one-to-one on each  $\omega^n, \,n \in \omega \}$.
\endroster

We order $S$ by $w \le v$ if dom$(v) \subseteq \text{dom} (w)$ and

\indent\indent (1) \quad $\langle \alpha, 0 \rangle \;\in \text{dom}(v)
\Rightarrow v(\alpha, 0) \subseteq w (\alpha, 0)$,

\indent\indent (2) \quad $\langle \alpha, 1 \rangle \;\in \text{dom}(v)
\Rightarrow
v(\alpha, 1) = w  (\alpha, 1)$,

\indent\indent (3)\quad
$\omega_2 \;\in \text{ dom }(v) \Rightarrow  v(\omega_2) \subseteq w
(\omega_2)$.

Let $G \subseteq S$ be a generic filter.  The forcing is designed to
add by finite conditions functions $k, h, g_{\alpha} (\alpha <
\omega_2)$: \newline
\noindent (1)\quad $g_{\alpha} = \bigcup \{w(\langle  \alpha, 0
\rangle): w\in G \} \in \omega^ \omega$,

\noindent (2)\quad $h : \omega_2 \to \omega$ defined by
$h(\alpha) = n$ if $ \langle \langle \alpha, 1 \rangle, n \rangle \;\in G$,

\noindent (3)\quad $k = \bigcup \{ w(\omega_2):w\in G\}$ is a function
from $\omega^{<\omega}$ which is one-to-one on each $\omega^n$.

\noindent Note that $g_{\alpha} : \alpha < \omega_2$ are mutually generic
Cohen reals.

Obviously, $B(S) = \C_{\omega_2}$ and it is enough to find a
complete embedding of $B(P) = B$ into $B(S)$.  We show
how to read off a $P$-generic sequence
$\langle f_{\alpha} : \alpha < \omega_2 \rangle$
from a generic filter $G \subseteq S$ and define the
associated projection $pr : B( S) \to B$
on a dense set $D \subseteq S$.  The
verification is left to the reader.

So let $G \subseteq S$ be generic, $k, h, g_{\alpha} \in V[G]$ as
above.  In $V [G]$, we define $f_{\alpha} : \omega \to \omega$ by $f_{\alpha}
\restriction h (\alpha) = g_{\alpha} \restriction h (\alpha), f_{\alpha} (n) =
k (g_{\alpha}
\upharpoonright (n+1))$ for $n \ge h (\alpha)$.

To define the projection we let $D$ be the set of all
$w \in S$ such that (1) $\langle \alpha, 0 \rangle  \;\in
\text{dom} (w) \text{ iff } \langle \alpha, 1 \rangle \;
\in \text{dom} (w)$, (2) the sequences  $w (\alpha, 0),
\;\langle  \alpha, 0 \rangle\in  \text{dom}(w)$,
are pairwise distinct and there is an  integer
$n \in \omega$ such that $w (\alpha, 0) \in \omega^n$,
and $\,w (\alpha, 1) < n$, and (3)
$\omega _2 \in \text{dom} (w) \text{ and dom}(w (\omega _2)) =
\{s \in \omega^{\le n} : s$ is an initial segment
of some  $w (\alpha, 0)\}$.

It is immediate that $D \subseteq S$ is dense.  We define $pr(w)$,
for  $w \in D$, to be the $z \in P$ for which:
\roster
\item  dom $(z) = \{\alpha: \;\langle \alpha, 0 \rangle \;\in
\text{ dom} (w)\}$

\item  $z (\alpha) \upharpoonright w (\alpha, 1) = w (\alpha, 0)
\upharpoonright w (\alpha, 1)$ 
\newline
and $z (\alpha) (m) =w (\omega_2)(w (\alpha, 0)) \upharpoonright (m+1~)
\text{ for } w (\alpha, 1) \le m < \text{ dom} (w(\alpha, 0))$.
$\square$
\endroster 
\enddemo

Semi-Cohen algebras with more complicated properties can be produced
using similar methods.  For example, the forcing from \cite{He} is similar
to the one we just described. The following argument describes
a semi-Cohen algebra $B$
of uniform density $\frak c^+$ that cannot be embedded into a Cohen algebra.

Let $P$ be a forcing for adding a modulo finite
increasing chain of functions from $\omega$ into the rationals $Q$ of 
length $\frak c^+$.  Then $B = B(P)$
is semi-Cohen and since no Cohen algebra adds such a chain \cite{Ku},
$B$ does not embed into any $\Ck$.  To be more precise, set
$P = \{z : \text{dom}(z) \in [\frak c^+]^{<\omega}$ and there is an
integer $n_z \in \omega$ such that for $\alpha \in \text{dom} (z)$
$z (\alpha) \in Q^{n_z}\}$.  The ordering is defined by $z \le w$ if
$z$ is a coordinatewise extension of $w$ and for $\alpha < \beta$ both
in dom$(w)$ and $n_w \le m < n_z$ we have $z (\alpha) (m) < z (\beta) (m)$.
Thus if $G \subseteq P$ is a generic filter and $\alpha < \frak c^+$, we can set
$f_{\alpha} : \omega \to Q$ to be $\bigcup \{f : \;\langle \alpha,
f \rangle \;\in G\}$ and
we will have $\alpha < \beta < \frak c^+ \Rightarrow f_{\alpha}(n) <
f_{\beta} (n)$ for all but finitely many $n \in \omega$.
The only thing
left to prove is the semi-Cohenness of $B = B(P)$.  For
$X \subseteq \frak c^+$ define $P_X = \{ z \in P : \text{dom}(z) \subseteq
X\}$.  Then $P_X$ is a regular subposet of $P$, with the projection
$pr: P \to P_X$ defined by $pr(z) = z\upharpoonright X$.  Thus
$\{P_X: X \in [\frak c^+]^{\omega}\}$ is a club set of regular
subposets of $P$ and consequently, $B=B(P)$ is semi-Cohen.

\medpagebreak

The last example in this section gives a sharper result and provides a
completely different semi-Cohen algebra.

\proclaim
{Theorem 5.5} \ There is a semi-Cohen algebra $B$ of uniform density
$\aleph_2$ which cannot be completely embedded into a Cohen algebra.
\endproclaim

\demo
{Proof} \ Let $P$ be Tennenbaum's forcing for adding a Souslin tree
with finite conditions \cite{T}.
So $P = \{ \langle t, <_t \rangle \;:
t \in [\omega_1]^{<\omega} \text{ and } <_t$
is a tree order on $t$ respecting the
ordering of ordinals$\}$;
we order $P$ by reverse extension.  Let the $P$-term for a Souslin tree
$\dot T$ be defined by $P \Vdash$ ``$\dot T =
\;\langle \omega_1, <_T \rangle$, where
$<_T = \bigcup \{ <_t : \langle t, <_t \rangle  \in G\}$''.\enddemo

\proclaim
{Lemma 5.6} \ (1) \ $B(P){\cong} \C_{\omega_1}$,

\hskip .7in (2) \ the completion of $P * (\dot T
\text{ upside down})$ is isomorphic to $\C_{\omega_1}$,

\hskip .7in (3) \ $B({P * \dot S}) \cong \C_{\omega_1}$,
where $\dot S$ is a $P$-name for the forcing 
that makes $\dot T$ special, cf. \cite{B-M-R}.
\endproclaim

\demo
{Proof} \ (1) is a consequence of both
(2) and (3).  The proofs of (2) and (3) are similar; we prove (3).

All we have to show by Theorem 4.2 is that the algebra $B({P * \dot S})$
is semi-Cohen.  The forcing $P * \dot S$ has a dense subset $D =
\{ \langle t, <_t, f \rangle \;: \; \langle t, <_t \rangle \;
\in P$ and $f : t \to \omega$ is a function such that
$\alpha <_t \beta \Rightarrow f(\alpha) \ne f (\beta) \}$,
ordered by reverse extension.  Let $D_{\alpha} =
\{ \langle t, <_t, f \rangle \; \in D : t \le \alpha\}$.

It is enough to show that for a limit ordinal $\alpha \in \omega_1$
the poset $D_{\alpha}$ is a regular subposet of $D$.  Then
$\{D_{\alpha} : \alpha \in \omega_1 \text{ limit}\}$ is a club
subset of regular subposets of $D$, proving the semi-Cohenness of
$B({P * \dot S})$.  So fix $\alpha \in \omega_1$ limit.  We
define the projection $pr: E \to D_{\alpha}$ for a dense set
$E \subseteq D.$ 
Let $E = \{ \langle t, <_t, f \rangle \; \in D :
\forall  \beta \in t \cap \alpha$ if
there is a $\gamma \in t$ such that $\beta <_t \gamma$ and $f(\gamma) = n$
then there is $\gamma \in t \cap \alpha$ such that
$\beta <_t \gamma$ and $f(\gamma) = n\}$.

Obviously, $E \subseteq D$ is dense and a function
$pr: E\to D_{\alpha}$ defined by
$pr \langle t, <_t, f \rangle \;= \ \langle t \cap \alpha,
\; <_t \upharpoonright \alpha^2,\ f \upharpoonright
\alpha \rangle$ is a projection.  $\square$
\enddemo

By a similar argument as in the Lemma, the algebra
$B({P * \dot X * \dot Q})$ is isomorphic to
$\C_{\omega_1}$, where
$\dot X$ is (forced to be) a four-element algebra with  two atoms $x,
1-x$ and $1 \Vdash_P x \Vdash_{\dot X}$
``$\dot Q$ is the forcing $\dot T$ upside
down'', $1 \Vdash_P 1-x \Vdash_{\dot X}$
``$\dot Q$ is the $\dot T$-specialization forcing.''

Finally, we are in a position to define the complete algebra $B$
witnessing the statement of the Theorem.  We let $B =B( R)$
where $R$ is the finite support iteration
$$
	R = \langle R_{\alpha} : \alpha \le \omega_2, \; \dot U_{\alpha}
: \alpha < \omega_2 \rangle
$$
where $U_0 = \C_{\omega_1}$, for $0 < \alpha < \omega_2$ we have
$B(R_{\alpha}) \cong \C_{\omega_1}$ and at stage $\alpha$,
$0 < \alpha < \omega_2$, we find an isomorphism $i_{\alpha} :
B( R_{\alpha}) \to B(P)$, get an $ R_{\alpha}$-
name $\dot T_{\alpha}$ for a Souslin tree and set
$\dot U_{\alpha} = X_{\alpha} * \dot Q_{\alpha}$ as above, that is,
$X_{\alpha}$ has two atoms $x_{\alpha}$ and $1 - x_{\alpha}$ and
$x_{\alpha}\Vdash_{X_{\alpha}}$ ``$\dot Q_{\alpha}$ is the forcing
with $\dot T_{\alpha}$ upside down'' and  $1 - x_{\alpha}
\Vdash_{X_{\alpha}}$ ``$\dot Q_{\alpha}$ is the
$\dot T_{\alpha}$-specialization forcing.''

We have $B = \bigcup_{\alpha < \omega_2} B_{\alpha}$, where
$B_{\alpha} = B(R_{\alpha})$.

\proclaim
{Lemma 5.7} \ $B$ is a semi-Cohen algebra of uniform density
$\aleph_2$.
\endproclaim

\demo
{Proof} \ By induction on $\alpha \in \omega_2$ we prove that
$B_{\alpha}$ is a semi-Cohen algebra of uniform density $\aleph_1$, i.e.
$B_{\alpha} \cong \C_{\omega_1}$.  At limits steps, we use Theorem 4.2,
and at successor steps, we apply the observation following the previous
Lemma.  Then Lemma follows immediately from Theorem 4.2. $\square$
\enddemo

\proclaim
{Lemma 5.8} \ $B$ cannot be completely embedded into a Cohen algebra.
\endproclaim

\demo
{Proof} \ Assume for contradiction that it can be.  Then $B$ can be
embedded into $\C_{\omega_2}$ and we can consider $B$ as a complete
subalgebra of $\C_{\omega_2}$.  The following can be proved by a
simple closure argument.
\enddemo

\proclaim
{Claim 5.9} \ The set $C = \{\alpha \in \omega_2 : B_{\alpha}
\subseteq \C_{\omega_2}\} \subseteq \omega_2$ is closed unbounded.
\endproclaim

Even better, we have

\proclaim
{Claim 5.10} \ The set $D = \{\alpha \in C$ : for no
condition $p \in C_{\alpha}$, $p$ decides the
statement ``$\dot x_{\alpha} \in G$''$\}$ contains a club
in $\omega_2$.
\endproclaim

\demo
{Proof} \ First, define a $B$-name $\dot f$ for a function from
$\omega_2$ to $\{0,1\}$ by $B \Vdash$ $``\dot f (\alpha) = 1$ iff
$\dot x_{\alpha} \in \dot G$''.  Obviously,
$$B \Vdash\text{``}\dot f\text { is a }\C_{\omega_2}\text{-generic
function over the ground
model''.}\tag {$*$}$$

Now assume that the set $C - D$ is stationary.  For any
$\alpha \in C - D$ let $p_{\alpha} \in C_{\alpha}$
decide the statement ``$x_{\alpha} \in G$''.
There exist a stationary set $S \subseteq
C - D$ and a $p \in \C_{\omega_2}$ such that for every
$\alpha \in S$ we have $p = p_{\alpha}$.  Then
$p \Vdash$ ``$\dot f \upharpoonright S$
belongs to the ground
model'', contradicting ($*$) since the set $S$ is infinite.
Hence $D$ contains a closed unbounded set.  $\square$
\enddemo

Now fix an ordinal $\alpha \in D$.  We have
$B_{\alpha} \le_{\text{reg}} \C_{\alpha}$ and so $\dot T_{\alpha}$
is a $\C_{\alpha}$-name for an $\omega_1$-tree.  We reach a
contradiction:

{\it Case I.}  There is $p \in C_{\alpha}$ such that
$p \Vdash_{\C_{\alpha}}$ ``$\dot T_{\alpha}$ has a
cofinal branch''.   Since $\alpha \in D$,  we have
$p \cdot (1- x_{\alpha}) \ne 0$ in $C_{\omega_2}$.  Also,
$1- x_{\alpha} \Vdash_B$ ``$\dot T_{\alpha}$ is a special
tree''.  By upwards absoluteness,
$p \cdot (1- x_{\alpha}) \Vdash_{\C_{\omega_2}}$
``$ T_{\alpha}$ is a special tree with a cofinal branch''.
It follows that $\omega_1$ must be collapsed, contradicting c.c.c.
of $\C_{\omega_2}$.

{\it Case II.}  $C_{\alpha} \Vdash$ ``$\dot T_{\alpha}$ has no
cofinal branches''.  Then  $\C_{\omega_2} \Vdash$
``$\dot T_{\alpha}$ has no cofinal branches'' since the residue
forcing $\C_{\omega_2} : \C_{\alpha}$ is Cohen and as such
does not  add branches to Aronszajn trees (Lemma 4.7).  However,
$x_{\alpha} \Vdash_B$ ``$\dot T_{\alpha}$ has a cofinal branch'',
so by upwards absoluteness  $x_{\alpha} \Vdash_{\C_{\omega_2}}$
``$\dot T_{\alpha}$ has a cofinal branch'', contradiction.

This completes the proof of Theorem 5.5. $\square$

\medpagebreak

The following example shows that Theorem 4.5 fails for
Cohen algebras:

\proclaim
{Theorem 5.11} \ There exists an increasing $\omega$-chain of the
Cohen algebras $B_n$, with $B_n\le_{\text{reg}} B_{n+1}$, whose
union is not a Cohen algebra.
\endproclaim

\demo
{Proof} \ Let $P$ be the forcing notion from Theorem 5.1.
A generic $G$ on $P$ yields $\omega_2$ functions
$\{g_{\alpha} : \alpha < \omega_2\}$ from $\omega$ into
$\omega$.  For each $n$ and each $\alpha < \omega_2$,
let $g_{\alpha,n} (k) = g_\alpha(k) \text{ mod } 2^n$, and let $G_n =
\{g_{\alpha, n} : \alpha < \omega_2\}$.

We have $V[G] = V [ \{G_n : n \in \omega \}]$, and for all $n<m$,
$V [G_n] \subset V [G_m]$.

Let $P_n$ be the forcing with finite conditions that adjoins
$\omega_2$ functions from $\omega$ into $2^n$.  $B (P_n)$ is a Cohen
algebra; we claim that $G_n$ is a generic on $P_n$:

If $D$ is an open dense set in $P_n$, let $E$ be the set of all
conditions whose projection belongs to $D$.  As $E$ is dense in
$P$, it follows that $G_n$ is generic on $P$.

Hence $B(P)$ is the limit of a regular $\omega$-chain of Cohen
algebras.  $\square$

The limit $B(P)$ in the above proof is embeddable into a Cohen
algebra.  By taking instead the forcing $P$ described in the
comments following the proof of Theorem 5.1, we obtain a regular
$\omega$-chain of Cohen algebras whose limit $B$ is not embeddable into
a Cohen algebra.  In this case, each $B_n$ has density
$\ge (2^{\aleph_0})^+$.
\enddemo

\head{6. Pseudo-Cohen algebras}\endhead

Looking for other classes of algebras which share some of the properties
of $\C_\kappa,$ we arrive at the following generalization of Definition 1.1:

\definition {Definition 6.1} An algebra $B$ of uniform density
is {\it pseudo-Cohen} if the set
$S=\{ A\in [B]^{\omega}:A\leq_{\text {reg}}B\}$ is stationary.
\enddefinition

By Lemma 3.4, pseudo-Cohen algebras are c.c.c.
and all reals added by them are in a Cohen-generic extension of the ground
model.
However, the class of pseudo-Cohen algebras does not have most of the closure
properties of the semi-Cohen class. While it is closed under regular
subalgebras (Theorem 6.2 below),
it is not closed under products or iterations,
since the nature of the witness set $S$ may vary.

\proclaim {Theorem 6.2} If $B$ is pseudo-Cohen and if $A$ is a regular
subalgebra of uniform density then $A$ is pseudo-Cohen.
\endproclaim

\demo {Proof} Follows closely the proof of Theorem 4.1. $\square$
\enddemo

The following generalizes Theorem 4.3:

\proclaim
{Theorem 6.3} The following are equivalent:

(a)\ $B$ is a pseudo-Cohen algebra.

(b)\ there is an $\aleph _0$-distributive
forcing $P$ such that $P\Vdash$``$B$ is Cohen".

(c)\ the first player does not have a winning strategy in the game $\Cal G.$
\endproclaim

\demo
{Proof} (a) implies (b). Let $B$ be a pseudo-Cohen algebra as witnessed by a
set $S$ and let $P$ be the standard forcing for shooting a club through $S,$
namely $P=\{ f:f$ is a function from some $\alpha+1<\omega _1$ to $S$ which
is increasing and continuous with respect to $\subset \},$ ordered by extension.
The forcing $P$ collapses $|B|$ to $\aleph _1$ and
is $\aleph _0$-distributive. In $V^P$ the algebra $B$ has size $\aleph_1$
and a club subset of regular subalgebras, therefore, it is Cohen by Theorem 4.2.

(b) implies (a). Let $P\Vdash$``$B$ is a Cohen algebra" for
some $\aleph_0$-distributive forcing $P.$ Fix a function $f:B^{<\omega}\to B.$
We must produce a countable regular subalgebra $A$ of $B$ closed under $f.$
Such an algebra certainly exists in the generic extension by $P,$ where $B$
is Cohen. But since $P$ does not add any countable subsets of $B,$ such a
subalgebra $A$ exists already in the ground model.

(a) implies (c). Let $\sigma$ be a strategy for the player I in
the game $\Cal G$
associated with a pseudo-Cohen algebra $B.$ Since the witness set $S$ is
stationary,
we can fix a large regular cardinal $\theta$ and a countable submodel
$M\prec\langle
H_\theta ,\in,B,\sigma\rangle$ with $M\cap B\leq_{\text{reg}}B.$
Since the model $M$ is closed under $\sigma,$ 
there is a play of $\Cal G$ in which player I uses
$\sigma$, and in which the second player picks all the elements of $M\cap B.$
This shows that $\sigma$ is not a winning strategy for the player I.

(c) implies (a): Left to the reader. 
\qed
\enddemo

Pseudo-Cohenness is preserved by proper forcing extensions:

\proclaim {Proposition 6.4} If $B$ is pseudo-Cohen and $W$ is a proper-forcing
extension of the universe then $W \Vdash B $ is pseudo-Cohen.
\endproclaim

\demo {Proof}
The witness $S$ from Definition 6.1 remains stationary in $W.$ $\square$
\enddemo

We finish this section by giving an example of a pseudo-Cohen algebra
that is not semi-Cohen. Fix a set
$\{l_\alpha:\alpha<\omega _1$ limit$\} ,$ where $l_\alpha:\omega\to\alpha$
is an increasing sequence of ordinals converging to $\alpha.$
Fix $S\subset \omega_1$ co-stationary and let $P_S=\{\langle f,s\rangle:f$
is a finite function from $\omega _1$ to $\{0,1\}$ and $s$ is a finite set of
countable ordinals$\}$ ordered by $\langle g,t\rangle\leq\langle f,s\rangle$
if $f\subset g,$ $s\subset t$ and for every $\beta\in\ $dom$(g- f)$
and every $\alpha\in s\cap S,$ if $\beta\in\ $ran$(l_\alpha)$ then $g(\beta
 )=0.$

In the generic extension by $P_S$, let
$F:\omega _1\to \{0,1\}$  be defined
by $F=\bigcup\{ f:\langle f,\emptyset\rangle$ is in the generic filter$\} .$
By the definition of the ordering on $P_S,$ a countable
limit ordinal $\alpha$ is in $S$ just in the case when only
finitely many $\beta$'s in ran$(l_\alpha)$ have $F(\beta )=1.$

We claim that $B(P_S)$ is a pseudo-Cohen algebra. To show
this, it is enough
to prove that for every $\alpha\in\omega_1- S$ the poset
$R_\alpha=\{\langle f,s\rangle\in P_S:s\subset\alpha,$
dom$(f)\subset\alpha\}\subset P_S$ is a regular subposet of $P_S.$
Fix an arbitrary condition $\langle f,s\rangle\in P_S.$ We shall produce
a condition $\langle g,t\rangle\in R_\alpha$ such that any extension of
it in $R_\alpha$ is compatible with $\langle f,s\rangle$ in $P_S.$
Let $x=s\cap S- \alpha$ and $y=\{ \beta <\alpha:\beta\notin$dom$(f)$
and $\exists \xi\in x\ \beta\in$ ran$(l_\xi)\}.$ Thus the set $y$ is finite,
because $\alpha\notin S.$ We let $g$ be the extension of
$f\restriction \alpha$ to $\alpha\cup y$ such that $g=0$ on $y$,
and $s=t\cap \alpha.$
It is easy to see that $\langle g,s\rangle$ is as required.

Thus $B(P_S)$ is a pseudo-Cohen algebra. If the set $S\subset \omega _1$ is
chosen stationary co-stationary, it is not difficult to prove that $B(P_S)$
is not semi-Cohen and $B(P_S)\times B(P_{\omega _1- S})$ is not
pseudo-Cohen. Thus we have:

\proclaim {Theorem 6.5} There exists a pseudo-Cohen algebra that is not
semi-Cohen.
\endproclaim

We can also show that the algebra constructed by Veli\v ckovi\'c
in [Ve] for adding a Kurepa tree by a c.c.c. forcing from $\square _{\omega_2}$
is also pseudo-Cohen.

\head{7.  Potentially Cohen algebras}\endhead

Generalizing properties (c) in Theorem 4.3 and (b) in Theorem 6.2 we arrive
at the following notion:

\definition
{Definition 7.1}  An algebra $B$ of uniform density
is {\it potentially Cohen} if
there is an $\omega_1$-preserving forcing notion such that
$Q \Vdash$ ``$B$ is a Cohen algebra''.
\enddefinition

Every pseudo-Cohen algebra is potentially
Cohen.  In this section we investigate the converse.
We prove:

\proclaim
{Theorem 7.2} \ (CH)  Every potentially Cohen algebra is
pseudo-Cohen.
\endproclaim

\proclaim
{Theorem 7.3} \ (MA+$\lnot$CH)  The measure algebra is potentially
Cohen.
\endproclaim

Notice that the measure algebra is not pseudo-Cohen,
since it does not add Cohen reals.  Thus Theorem 7.3 proves the
necessity of the CH assumption in Theorem 7.2 and shows that in
general, potential Cohenness of an algebra is a considerably weaker
property. Indeed, the only significant properties
of an algebra $B$ implied by its potential Cohenness
without further assumptions seem to be included in the
following simple lemma.

\proclaim
{Lemma 7.4} \ If $B$ is a  potentially Cohen algebra, then $B$
is c.c.c. and not $\aleph_0$-distributive.
\endproclaim

\demo
{Proof} \ Let $B$ be potentially Cohen and $Q \Vdash$
``$B$ is Cohen'', for some $\omega_1$-preserving
forcing $Q$.  If $A \subseteq B$ is an uncountable antichain,
then $A$ remains uncountable in $V^Q$, contradicting c.c.c. of $B.$
So the algebra $B$ is c.c.c.  If $B$ were
$\aleph_0$-distributive, then there would be a Souslin tree $T$
such that $B(T) \le_{\text{reg}} \overline B$.  Then
$Q \Vdash$ ``$B(T)\le_{\text{reg}} \overline B \cong
\C_\kappa$ for some $\kappa$'' and since the uniform density of
$B(T)$ is $\aleph_1$, by Theorem 4.2  $Q \Vdash B(T) \cong
\C_{\omega_1}$.  So in the generic extension
by $Q$, forcing with the $\omega_1$-tree $T$ adds reals.
Consequently $Q \Vdash$ ``$B(T)$ collapses $\omega_1$ and
so $B(T) \not\cong \C_{\omega_1}$'' contradiction.
$\square$
\enddemo

\demo
{Proof of Theorem 7.2} \ Assume that $B$ is a potentially Cohen
poset as witnessed by an $\omega_1$-preserving forcing $Q$.
First, we treat the special case when $|B| = \aleph_1$.  Let us
enumerate $B = \{p_{\alpha} : \alpha < \omega_1\}$.  We shall show
that the set $S =
\{\beta < \omega_1 :\ B_{\beta} =
\{p_{\alpha} :$ $\alpha < \beta\}\ \le_{\text{reg}} B\}$ is a
stationary subset of $\omega_1$, which proves the pseudo-Cohenness
of the poset $B$.  Indeed, in the forcing extension $V^Q$ we will
have $S = \{\beta < \omega_1 : B_{\beta} =
\{p_{\alpha} : \alpha < \beta\}\le_{\text{reg}} B\}$ and since
$V^Q \vDash$ ``$ B$ is a Cohen algebra'', the set $S$ contains a
closed unbounded subset in $V^Q$.  Consequently, $S$ is stationary in $V$.

In the case of a poset $B$ of higher cardinality, we have to
prove that the set $S = \{A \in [B]^{\omega} :
A \le_{\text{reg}} B\}$ is stationary.  So let $f : B^{<\omega}
\to B$ be an arbitrary function.  We shall produce an element
of the set $S$ closed under $f$.
\enddemo

\proclaim
{Lemma 7.5} \ (CH)  There is a regular subalgebra $A_1$ of $B$ of
size $\aleph_1$ closed under the function $f$.
\endproclaim

\demo
{Proof} \ Choose a large regular cardinal $\kappa$ and a submodel
$M \prec H_{\kappa}$ of size $\aleph_1$ such that $f$, $B \in M$
and $[M]^{\omega} \subset M$.  This is possible since the
Continuum Hypothesis holds.  Now any maximal antichain $X$ of the
poset $B \cap M$ is an antichain in $B$ and so is countable.
Therefore, $X \in M$ and by elementarity $X$ is a maximal
antichain of $B$.  Consequently, the algebra $A_1 = M \cap
B \le_{\text{reg}} B$ is as required. $\square$

Let $A_1$ be as in Lemma 7.5.  Now
$Q \Vdash$ ``$A_1\le_{\text{reg}} B \le_{\text{reg}}\overline B$
and so $A_1$ is a regular subalgebra of a Cohen algebra and by
Theorem 4.2, $\overline A_1$ is a Cohen algebra itself''.  Thus
$A_1$ is potentially Cohen.  By the first part of the proof, $A_1$
is pseudo-Cohen and there is a countable subposet
$A\le_{\text{reg}} A_1$ closed under the function $f$.  So,
$A \le_{\text{reg}} A_1 \le_{\text{reg}} B$ is as required. $\qed$
\enddemo

\demo
{Proof of Theorem 7.3} \ Assume $MA + \frak c > \aleph_1$.  The theorem will
follow from these two lemmas:
\enddemo

\proclaim
{Lemma 7.6} \ ($MA + \neg CH$)  There is an $\omega_1$-preserving
forcing $Q$ such that
$Q \Vdash$``{\rm cf}$(\frak c^V) = \omega$''.
\endproclaim

\proclaim
{Lemma 7.7} \ ($MA + \neg CH$)  For any family
$\{B_{\alpha} : \alpha < \kappa\}$ of positive Borel sets,
where $\kappa < \frak c$, there are positive closed sets
$C_i : i < \omega$ such that $\forall \alpha < \kappa \,
\exists \, i < \omega \,\,C_i \subseteq B_{\alpha}$.
\endproclaim

Let a forcing $Q$ be as in Lemma 7.6.  Then
$Q \Vdash$ ``the measure algebra $B$ from the ground model
has a countable dense set''.  To see this, in the ground model we
enumerate $B$ as $\langle B_{\alpha} : \alpha < \frak c \rangle$,
choosing Borel representatives for each equivalence class.  By
Lemma 7.7, for each ordinal $\kappa < \frak c$ there are closed positive
sets $C^{\kappa}_i : i < \omega$ such that
$\forall \alpha < \kappa \,\exists \, i < \omega \,\, C^{\kappa}_i
\subseteq B_{\alpha}$.  Now if $Q \Vdash$
``$\langle \kappa_j : j < \omega \rangle$
is a sequence converging to $\frak c^V$'', then
$Q \Vdash$ ``$\{ C^{\kappa_j}_i : i, j < \omega \}$
is a countable dense subset of $B$''.  Thus the forcing $Q$ witnesses
the fact that the measure algebra $B$ is potentially Cohen.

\demo
{Proof of Lemma 7.6} \ Let
$\frak I = \{ A \subset \Cal P_{\aleph_2}
(\frak c) : \exists y \in \Cal P_{\aleph_2} (\frak c) \;\forall
x \in A \,\, y \not\subset x\}$ be the ideal of bounded subsets of
$\Cal P_{\aleph_2} (\frak c)$.  Since $\frak c > \aleph_1$, the
ideal $\frak I$ is proper and $\aleph_2$-complete.  We define the
forcing $Q$ to be the set of all trees $T$ of finite sequences of
elements of $\Cal P_{\aleph_2} (\frak c)$ such that the tree $T$
has a trunk $t$ and for each sequence $s \in T$ extending $t$ the
set $A_s = \{y \in \Cal P_{\aleph_2}  (\frak c) : s^{\frown}
\langle y \rangle \in T\}$ is not in the ideal $\frak I$.  In the
spirit of Namba forcing proofs, one can argue that $Q$ preserves
$\aleph_1$.  Also, if $G \subset Q$ is a generic filter, then in
$V [G]$, the set $\bigcup G$ is an $\omega$-sequence of sets which are
of cardinality $\aleph_1$ in the ground model and whose union
exhausts all of $\frak c^V$.  Since by $MA + \neg CH$, cf$(\frak c) >
\omega_1$ in $V$, the suprema of these sets are smaller than
$\frak c^V$ and converge to $\frak c^V$.  $\qed$
\enddemo

\demo
{Proof of Lemma 7.7} \ Let us recall one of the
definitions of the amoeba forcing:
$$A=\{\langle \Cal O, \delta \rangle : 0 < \delta \le 1 \text{ is a
real number and }\Cal O \subset [0, 1] \text{ is an open set of measure}
< \delta\}.$$
The order is defined by $\langle \Cal O, \delta
\rangle \ge \langle \Cal P, \gamma \rangle$ if $\gamma \le \delta$
and $\Cal O \subset \Cal P$.  The forcing $A$ is known to be
$\sigma$-linked and so a finite support product $A^{\omega}$
of $\omega$ copies of $A$ is c.c.c.

Fix a family $\{B_{\alpha} : \alpha < \kappa \}$ of positive
Borel sets of reals, where $\kappa < \frak c$.  For $\alpha <
\kappa$ we define sets $D_{\alpha} \subset A^{\omega}$ by
$p \in D_{\alpha}$ iff $\exists \,i < \omega \,p(i) =
\langle \Cal O, \delta \rangle$ for some real $\delta \in
(0, 1]$ and an open set $\Cal O$ such that $[0, 1] -
\Cal O \subset B_{\alpha}$.  It is easy to see that for $\alpha <
\kappa$ the set $D - \alpha \subset A^{\omega}$ is open dense:
if $p \in A^{\omega}$ and $\alpha < \kappa$, then one can
choose an integer $i$ with $i \not\in$ support$(p)$ and a closed
positive set $C \subseteq B_{\alpha}$.  Then $q = p \cup
\{ \langle i, \langle [0, 1] - C, 1 \rangle \rangle \}$
is a condition in $D_{\alpha}$ which is smaller than $p$.

Also, for an integer $i \in \omega$ define a dense subset $E_i
\subset A^{\omega}$ by $p \in E_i$ if $p(i) = \langle \Cal O,
\delta \rangle$ for some $\delta < 1$.

By Martin's Axiom there is a filter $G \subset A^{\omega}$
meeting all the dense sets $D_{\alpha} : \alpha < \kappa$ and
$E_i : i < \omega$.  For $i < \omega$ we define a closed set
$C_i = [0, 1] - \bigcup \{ \Cal O : \langle i, \langle
\Cal O, 1 \rangle \rangle \in G\}$.  Since the filter $G$ meets
all the $E_i$'s, the sets $C_i$ are positive.  Since the filter
meets all the sets $D_{\alpha}$, for every $\alpha < \kappa$ there
is an integer $i$ such that $C_i \subset B_{\alpha}$.  $\qed$
\enddemo

\enddocument
\end